\documentclass[11pt]{article}
\usepackage[francais]{babel}
\usepackage[latin1]{inputenc}
\usepackage[T1]{fontenc}
\usepackage{a4,amsmath,amsfonts}
\usepackage{amssymb}
\title{Théorème de la base normale effectif}
\author{Pascal Autissier}
\begin{document}

\maketitle

\newcommand{\D}{\displaystyle}
\newcommand{\sm}{\smallsetminus}

\textbf{English title:} Effective normal basis theorem.\\

\textbf{Abstract:} Let $K$ be a finite Galois extension of $\mathbb{Q}$. The
normal basis theorem provides an element of $K$ whose conjugates form a
$\mathbb{Q}$-basis of $K$. Here we obtain such an element with controlled size.
This improves a recent result by Fukshansky and Jeong. By the way, we estimate
Minkowski's minima of ideals of integers of number fields.\\

\textbf{Résumé:} Soit $K$ une extension finie galoisienne de $\mathbb{Q}$. Le
théorème de la base normale fournit un élément de $K$ dont les conjugués forment
une $\mathbb{Q}$-base de $K$. On obtient ici un tel élément de taille contrôlée.
Cela améliore un résultat récent de Fukshansky et Jeong. On estime au passage
les minima de Minkowski des idéaux d'entiers des corps de nombres.\\


\section{Introduction}

Soit $K$ une extension finie galoisienne de $\mathbb{Q}$, de groupe de Galois
$G=\{\sigma_1,\dots,\sigma_n\}$. Le théorème de la base normale (voir par
exemple le théorème 28 de \cite{Art}) énonce l'existence d'un $\alpha\in K$ tel
que $(\sigma_1(\alpha),\dots,\sigma_n(\alpha))$ soit une $\mathbb{Q}$-base de
$K$. On montre dans le présent article comment contrôler la ``taille'' d'une
telle base en fonction du degré $n$ et du discriminant $D_K$ de $K$.\\

\textbf{Théorème 1:} \textit{Il existe $\alpha\in O_K$ tel que
$(\sigma_1(\alpha),\dots,\sigma_n(\alpha))$ soit une $\mathbb{Q}$-base de $K$ et
que $|\sigma_i(\alpha)|\leqslant n|D_K|^{1/n}$ pour tout $i\in[[1,n]]$.}\\

En particulier, la hauteur de cet élément vérifie
\[h(\alpha)=\frac{1}{n}\sum_{i=1}^n\ln\max\bigl(1,|\sigma_i(\alpha)|\bigr)\leqslant\frac{\ln|D_K|}{n}+\ln n\ .\]
\`A titre de comparaison, Fukshansky et Jeong \cite{FuJe} obtiennent une base
$(\sigma_1(\alpha'),\dots,\sigma_n(\alpha'))$ avec $\alpha'\in K$ et
$h(\alpha')\leqslant(n-1)(4n-3)\ln|D_K|+ c(n)$, où $c(n)\sim12n\ln n$.\\

On observe également qu'une base normale ne peut pas être trop petite:\\

\textbf{Proposition 2:} \textit{Soit $\beta\in O_K$ tel que
$(\sigma_1(\beta),\dots,\sigma_n(\beta))$ forme une $\mathbb{Q}$-base de $K$.
Alors $|\sigma_1(\beta)|^2+\cdots+|\sigma_n(\beta)|^2\geqslant|D_K|^{1/n}$.}\\

L'argument du théorème 1 repose sur le fait que l'ensemble des $\alpha\in K$
tels que $(\sigma_1(\alpha),\dots,\sigma_n(\alpha))$ soit liée, est contenu dans
une hypersurface du $\mathbb{Q}$-espace $K$. On se ramène ainsi à trouver un
petit point du réseau $O_K$ qui évite cette hypersurface.\\

Après quelques préliminaires sur les réseaux (section 2), on estime à la section
3 les minima de Minkowski des idéaux d'entiers des corps de nombres. On en
déduit les résultats ci-dessus à la section 4.

Pour finir, on montre que notre approche permet aussi de donner une version
effective du théorème de l'élément primitif (voir la section 5).

\section{Géométrie des nombres}

\textbf{Définition:} Soit $E$ un $\mathbb{Q}$-espace vectoriel de dimension
finie $n\geqslant1$. Une application $f:E\rightarrow\mathbb{C}$ est dite
\textbf{polynomiale} de degré $d\in\mathbb{N}$ lorsque pour toute
base $(e_1,\dots,e_n)$ de $E$, il existe $P\in\mathbb{C}[X_1,\dots,X_n]$ de
degré $d$ tel que $f(a_1e_1+\cdots+a_ne_n)=P(a_1,\dots,a_n)$ pour tout
$(a_1,\dots,a_n)\in\mathbb{Q}^n$.\\

Soient $V$ un $\mathbb{R}$-espace vectoriel de dimension finie $n\geqslant1$ et
$\|\ \|$ une norme sur $V$. Lorsque $\Gamma$ est un réseau de $V$ et
$i\in[[1,n]]$, on note $\lambda_i(\Gamma)$ le $i$-ème minimum de Minkowski de
$\Gamma$, \textit{i.e.} $\lambda_i(\Gamma)=\min\bigl\{\max(\|e_1\|,\dots,\|e_i\|)\ ;\ (e_1,\dots,e_i)\textrm{ est une famille libre de }\Gamma\bigr\}$.\\

On aura besoin du lemme d'évitement suivant.\\

\textbf{Proposition 3 (Gaudron-Rémond):} \textit{Soit $\Gamma$ un réseau de $V$.
Soit $f:\Gamma_\mathbb{Q}\rightarrow\mathbb{C}$ une application non nulle,
polynomiale de degré $d$. Il existe $\alpha\in\Gamma$ tel que $f(\alpha)\neq0$
et que $\|\alpha\|\leqslant d\lambda_n(\Gamma)$.}\\

\textit{Démonstration:} C'est un cas particulier du théorème 1.1 de \cite{GaRe}.
Pour le confort du lecteur, on donne un argument direct. Par définition de
$\lambda_n(\Gamma)$, le réseau $\Gamma$ contient une famille libre
$(e_1,\dots,e_n)$ telle que $\|e_i\|\leqslant\lambda_n(\Gamma)$ pour tout
$i\in[[1,n]]$. Il existe un polynôme $P\in\mathbb{C}[X_1,\dots,X_n]$ de degré
$d$ vérifiant $f(a_1e_1+\cdots+a_ne_n)=P(a_1,\dots,a_n)$ pour tout
$(a_1,\dots,a_n)\in\mathbb{Q}^n$.

Le théorème des zéros combinatoire (théorème 1.2 de \cite{Alo}) fournit un
$(a_1,\dots,a_n)\in\mathbb{N}^n$ tel que $P(a_1,\dots,a_n)\neq0$ et que
$a_1+\cdots+a_n\leqslant d$. Alors l'élément $\alpha=a_1e_1+\cdots+a_ne_n$
convient. $\square$

\section{Minima des idéaux d'entiers}

\textbf{Lemme:} \textit{Soient $K$ un corps et $L$ une extension finie de $K$
de degré $n$. Soient $(x_1,\dots,x_k)$ et $(y_1,\dots,y_\ell)$ deux familles
$K$-libres de $L$. Supposons $k+\ell\geqslant n+1$. Alors la partie
$\bigl\{x_iy_j\ ;\ i\in[[1,k]]\textrm{ et }j\in[[1,\ell]]\bigr\}$ engendre le
$K$-espace vectoriel $L$.}\\

\textit{Démonstration:} \'Etant donnée une forme $K$-linéaire
$\varphi:L\rightarrow K$ non nulle, on va prouver qu'il existe $i\in[[1,k]]$ et
$j\in[[1,\ell]]$ tels que $\varphi(x_iy_j)\neq0$.

L'application $K$-linéaire $f:\begin{array}{ccc}
L&\rightarrow&L^\vee\\
x&\mapsto&\bigl(y\mapsto\varphi(xy)\bigr)\\
\end{array}$ est injective, donc la famille
$(f(x_1),\dots,f(x_k))$ est $K$-libre dans $L^\vee$. Le sous-espace
$\D V=\bigcap_{i=1}^k\mathrm{Ker}f(x_i)$ est ainsi de dimension
$n-k\leqslant\ell-1$. D'où l'existence d'un indice $j\in[[1,\ell]]$ vérifiant
$y_j\notin V$, puis de $i\in[[1,k]]$ tel que $\varphi(x_iy_j)=f(x_i)(y_j)\neq0$.
$\square$\\

Soit $K$ un corps de nombres de degré $n$. On note $r_1$ le nombre de
plongements réels de $K$ et $2r_2$ le nombre de ses plongements imaginaires.
Numérotons $\sigma_1,\dots,\sigma_n$ les plongements complexes de $K$ de sorte
que $\sigma_i(K)\subset\mathbb{R}$ pour $i\in[[1,r_1]]$ et que
$\sigma_{j+r_2}=\overline{\sigma_j}$ pour $j\in[[r_1+1,r_1+r_2]]$.

On munit $V=\mathbb{R}^{r_1}\times\mathbb{C}^{r_2}$ de la mesure de Lebesgue et on
plonge $K$ dans $V$ via $(\sigma_1,\dots,\sigma_{r_1+r_2})$. Si $I\subset K$ est
un idéal fractionnaire de $O_K$ de norme $N(I)$, on sait que $I$ est un réseau
de $V$ de covolume $\D\frac{\sqrt{|D_K|}}{2^{r_2}}N(I)$.

Soit $\|\ \|$ une norme sur $V$ telle que $\|xy\|\leqslant\|x\|\|y\|$ pour tout
$(x,y)\in V^2$. Désignons par $B$ la boule unité fermée de $V$.\\

\textbf{Proposition:} \textit{Soit $(k,\ell)\in[[1,n]]^2$ tel que
$k+\ell\geqslant n+1$. Soient $I$ et $J$ deux idéaux fractionnaires de $O_K$.
Alors
$\lambda_n(IJ)\leqslant\lambda_k(I)\lambda_\ell(J)$.}\\

\textit{Démonstration:} Il existe une famille $\mathbb{Q}$-libre
$(x_1,\dots,x_k)$ de $I$ vérifiant $\|x_i\|\leqslant\lambda_k(I)$ pour tout
$i\in[[1,k]]$. De même, il existe une famille libre
$(y_1,\dots,y_\ell)\in J^\ell$ telle que $\|y_j\|\leqslant\lambda_\ell(J)$ pour
tout $j\in[[1,\ell]]$.

La partie $\bigl\{x_iy_j\ ;\ (i,j)\in[[1,k]]\times[[1,\ell]]\bigr\}$
engendre le $\mathbb{Q}$-espace vectoriel $K$ d'après le lemme précédent. Et
pour tout $(i,j)\in[[1,k]]\times[[1,\ell]]$, on a
\[\|x_iy_j\|\leqslant\|x_i\|\|y_j\|\leqslant\lambda_k(I)\lambda_\ell(J)\ .\]
D'où l'inégalité énoncée. $\square$\\

\textbf{Corollaire:} \textit{Soit $I$ un idéal fractionnaire de $O_K$. On a la
majoration
\[\lambda_n(I)^n\leqslant\frac{4^{r_1+r_2}}{\mathrm{vol}(B)^2}|D_K|N(I)\ .\]}\\

\textit{Démonstration:} La proposition précédente donne
$\lambda_n(I)\leqslant\lambda_i(I)\lambda_{n+1-i}(O_K)$ pour tout $i\in[[1,n]]$,
donc
\[\lambda_n(I)^n\leqslant\prod_{i=1}^n\lambda_i(I)\prod_{j=1}^n\lambda_j(O_K)\ .\]
Appliquons le deuxième théorème de Minkowski à $I$ et à $O_K$:
\[\prod_{i=1}^n\lambda_i(I)\leqslant2^n\frac{\mathrm{covol}(I)}{\mathrm{vol}(B)}=\frac{2^{r_1+r_2}}{\mathrm{vol}(B)}\sqrt{|D_K|}N(I)\quad\textrm{et}\quad\prod_{j=1}^n\lambda_j(O_K)\leqslant\frac{2^{r_1+r_2}}{\mathrm{vol}(B)}\sqrt{|D_K|}\ .\]
Le résultat en découle. $\square$\\

\textbf{Exemple 4:} Notons $\|\ \|_\infty$ la norme sur $V$ définie par
$\|x\|_\infty=\max(|x_1|,\dots,|x_{r_1+r_2}|)$ pour
$x=(x_1,\dots,x_{r_1+r_2})\in\mathbb{R}^{r_1}\times\mathbb{C}^{r_2}$. En prenant
$\|\ \|=\|\ \|_\infty$, le volume de $B$ vaut $2^{r_1}\pi^{r_2}$; on obtient
$\D\lambda_n(I)^n\leqslant\Bigl(\frac{2}{\pi}\Bigr)^{2r_2}|D_K|N(I)$.

On retrouve ainsi une variante de la proposition 4.2 de \cite{Bay} via une
approche différente (voir aussi le théorème 1.6 de \cite{BS3TZ} pour une version
moins précise).

\section{Base normale}

Soit $K$ une extension finie galoisienne de $\mathbb{Q}$, de groupe de Galois
$G=\{\sigma_1,\dots,\sigma_n\}$. On désigne par
$\mathrm{Tr}:K\rightarrow\mathbb{Q}$ la forme $\mathbb{Q}$-linéaire trace.
Choisissons une $\mathbb{Z}$-base $(e_1,\dots,e_n)$ de $O_K$. Pour tout
$x\in K$, on pose
$\Delta(x)=\det\bigl[\mathrm{Tr}(e_i\sigma_j(x))\bigr]\in\mathbb{Q}$.\\

\textbf{Lemme:} \textit{Soit $x\in K$. On a $\Delta(x)\neq0$ si et seulement si
$(\sigma_1(x),\dots,\sigma_n(x))$ est une $\mathbb{Q}$-base de $K$.}\\

\textit{Démonstration:} Il suffit de remarquer que la forme
$\mathbb{Q}$-bilinéaire $(x,y)\mapsto\mathrm{Tr}(xy)$ est non dégénérée sur $K$.
$\square$\\

Prouvons maintenant les résultats de l'introduction.\\

\textit{Démonstration du théorème 1:} L'application $\Delta$ est non nulle
d'après le théorème de la base normale, et polynomiale de degré $n$: si
$(e'_1,\dots,e'_n)$ est une $\mathbb{Q}$-base de $K$, alors
$\Delta(a_1e'_1+\cdots+a_ne'_n)=P(a_1,\dots,a_n)$ pour tout
$(a_1,\dots,a_n)\in\mathbb{Q}^n$, avec
\[P=\det\biggl[\sum_{k=1}^n\mathrm{Tr}(e_i\sigma_j(e'_k))X_k\biggr]\in\mathbb{Q}[X_1,\dots,X_n]\ .\]
La proposition 3 fournit un $\alpha\in O_K$ tel que $\Delta(\alpha)\neq0$ et que
$\|\alpha\|_\infty\leqslant n\lambda_n(O_K)$. On conclut via la majoration de
l'exemple 4. $\square$\\

\textit{Démonstration de la proposition 2:} En posant
$\Gamma=\mathbb{Z}\sigma_1(\beta)+\cdots+\mathbb{Z}\sigma_n(\beta)$, on a la
relation $\bigl|\det\bigl[\sigma_i(e_j)\bigr]\bigr|\#(O_K/\Gamma)=\bigl|\det\bigl[\sigma_i\sigma_j(\beta)\bigr]\bigr|$, donc
\[\sqrt{|D_K|}\leqslant\bigl|\det\bigl[\sigma_i\sigma_j(\beta)\bigr]\bigr|\leqslant\bigl(|\sigma_1(\beta)|^2+\cdots+|\sigma_n(\beta)|^2\bigr)^{n/2}\]
par l'inégalité d'Hadamard. $\square$

\section{\'Elément primitif}

Soit $K$ un corps de nombres. Notons $\sigma_1,\dots,\sigma_n$ les plongements
complexes de $K$. Le théorème de l'élément primitif assure l'existence d'un
$\alpha\in K$ tel que $K=\mathbb{Q}(\alpha)$. Voici une version effective de cet
énoncé:\\

\textbf{Proposition 5:} \textit{Il existe $\alpha\in O_K$ vérifiant
$K=\mathbb{Q}(\alpha)$ et $|\sigma_i(\alpha)|\leqslant(n-1)|D_K|^{1/n}$ pour tout
$i\in[[1,n]]$.}\\

\textit{Démonstration:} Posons
$\D g(x)=\prod_{i=2}^n(\sigma_1(x)-\sigma_i(x))$ pour tout $x\in K$; observons
que $g(x)\neq0$ si et seulement si $K=\mathbb{Q}(x)$. L'application
$g:K\rightarrow\mathbb{C}$ est non nulle grâce au théorème de l'élément
primitif, et polynomiale de degré $n-1$. On conclut comme précédemment,
\textit{i.e.} en combinant la proposition 3 et l'exemple 4. $\square$

\ \\

Pascal Autissier. I.M.B., université de Bordeaux, 351, cours de la
Libération, 33405 Talence cedex, France.

pascal.autissier@math.u-bordeaux.fr

\end{document}